\documentclass{amsart}

\usepackage{amssymb,amsfonts}

\usepackage[dvips]{graphicx}

\theoremstyle{plain}

\newtheorem{thm}{Theorem}[section]
\newtheorem*{thm_sixtheoremsharp}{Theorem \ref{sixtheoremsharp}}
\newtheorem*{thm_uniqueness}{Theorem \ref{uniqueness}}
\newtheorem*{thm_widthlessthantwo}{Theorem \ref{widthlessthantwo}}
\newtheorem*{thm_wge1}{Theorem \ref{wge1}}
\newtheorem*{thm_knotsfan}{Theorem \ref{knotsfan}}
\newtheorem*{thm_we1means3gon}{Theorem \ref{we1means3gon}}

\newtheorem{cor}[thm]{Corollary}
\newtheorem*{cor_uniquenesscor}{Corollary \ref{uniquenesscor}}
\newtheorem*{cor_knotsfancor}{Corollary \ref{knotsfancor}}

\newtheorem{lem}[thm]{Lemma}

\newtheorem{example}[thm]{Example}

\theoremstyle{definition}
\newtheorem{defn}[thm]{Definition}

\theoremstyle{remark}

\newcommand{\Sth}{\textbf{S}$^3$}

\newcommand{\Hth}{\textbf{H}$^3$}

\begin{document}


\title[Totally Geodesic Seifert Surfaces in Knot and Link Complements II] {Totally Geodesic Seifert Surfaces in Hyperbolic Knot and Link Complements II}

\date{\today}

\author[Adams]{Colin Adams}
\author[Bennett]{Hanna Bennett}
\author[Davis]{Christopher Davis}
\author[Jennings]{Michael Jennings}
\author[Novak]{Jennifer Novak}
\author[Perry]{Nicholas Perry}
\author[Schoenfeld]{Eric Schoenfeld}

\address{Colin Adams, Department of Mathematics and Statistics, Williams College, Williamstown, MA 01267}
\email{Colin.C.Adams@williams.edu}
\address{Hanna Bennett, Department of Mathematics, University of Chicago,5734 S. University Avenue, Chicago, Illinois 60637}
\email{hbennett@uchicago.edu}
\address{Christopher Davis, Department of Mathematics, Massachusetts Institute of Technology, 77 Massachusetts Avenue, Cambridge, MA 02139}
\email{davis@math.mit.edu}
\address{Michael Jennings, PO Box 144, Ithaca, NY 14851}\email{mvj3@cornell.edu}
\address{Nicholas Perry, 50 Taft Avenue, Newton, MA 02465}\email{nperry@wso.williams.edu}
\address{Jennifer Novak and Eric Schoenfeld, Department of Mathematics, Stanford University, 450 Serra Mall, Bldg. 380, Stanford, CA 94305-2125}
\email{novak@math.stanford.edu}\email{erics@math.stanford.edu}

\begin{abstract}
We generalize the results of \cite{AS}, finding large classes of totally
geodesic Seifert surfaces in hyperbolic knot and link complements, 
each the lift of a rigid $2$-orbifold embedded in some hyperbolic $3$-orbifold.  In addition,
we provide a uniqueness theorem and demonstrate that many knots 
cannot possess totally geodesic Seifert surfaces by giving bounds on the width invariant
in the presence of such a surface.  Finally, we utilize these examples
to demonstrate that the Six Theorem is sharp for knot complements in 
the 3-sphere.
\end{abstract}

\maketitle

\footnotetext[1]{2000 Mathematics Subject Classification 57M25; 57M50}

\section{Introduction} \label{S:intro}

Define a knot or link in $S^3$ to be hyperbolic if its complement is a
hyperbolic 3-manifold. This implies that there is a covering map p from
{\Hth} to $S^3-K$ such that the covering translations are isometries of
{\Hth}. We say that an embedded or immersed surface $S$ in $S^3-K$ is
totally geodesic if it is isotopic to a surface that lifts to a set of
geodesic planes in {\Hth}.  Throughout this paper, we will be using the
upper half-space model of {\Hth}, where the lift of a particular
cusp neighborhood is a union of horoballs.  In particular, we employ
pictures generated by Jeff Weeks' program SnapPea \cite{Weeks} displaying
the pattern of horoballs in the cusp lifts by looking down at the
$\{xy\}$-plane from above.

In \cite{AS}, the first examples of totally geodesic Seifert surfaces in
knot complements were produced.  These examples were generated using rigid
$2$-orbifolds embedded in hyperbolic $3$-orbifolds.  The main idea is that
certain knot complements cover hyperbolic $3$-orbifolds -- if a surface
$S$ in the knot complement projects to a rigid 2-orbifold under the
covering map, then $S$ must indeed be totally geodesic.

Section \ref{S:generating} generalizes this class of examples,
utilizing spherical 3-orbifolds as listed in  \cite{Dunbar}, and rigid $2$-orbifolds 
as appear in  \cite{Thurston} .  A spherical $3$-orbifold $O$ has universal cover $S^3$, so if we drill out appropriate curves from $O$ to get a hyperbolic $3$-orbifold $O'$, then $O'$ will lift
to a hyperbolic knot or link complement in $S^3$.  Any rigid
$2$-orbifold embedded in $O'$ will lift to a totally geodesic surface in
the knot or link complement. We will be particularly interested in
surfaces, the boundary of which is the knot or link.

\begin{defn}
A \textbf{Seifert surface} in a knot or link complement is an orientable
surface whose boundary is the knot or link. A nonorientable surface with
boundary the knot or link is called a \textbf{nonorientable Seifert
surface}.
\end{defn}

We often consider Seifert surfaces in the knot or link exterior,
$S^3 - N(L)$. In this case, the boundary of $S$ is a union of
\textbf{$l$-curves}, one in each cusp boundary, where an
$l$-curve is defined to be a closed curve in the boundary of a cusp 
neighborhood
which intersects the meridian exactly once.  It is a fact that the
boundary of a Seifert surface in a knot complement is a
longitude, defined as the $l$-curve which has linking number $0$ with
the missing core curve of the cusp.

Additionally, we often need to distinguish between different types of
surfaces using the following categorization:

\begin{defn}
Let $S$ be an embedded surface in the complement of a link $L$. Then $S$ is
  \textbf{free} if $S^3 - N(L) - N(S)$ is a handlebody.
We say $S$ is \textbf{totally knotted} if $S^3 - N(L) - N(S)$ has
incompressible boundary. We say $S$ is \textbf{semifree} if there exists
a compressing disk for $\partial(S^3 - N(L) - N(S))$.  Note that
free implies semifree.
\end{defn}

In Subsection \ref{SS:examples} we provide a series of examples of
orientable and non-orientable Seifert surfaces (both free and totally
knotted) in knot and link complements using the methods described
in Section \ref{S:generating}.  It remains an open question as to
whether a knot can have a non-orientable totally geodesic Seifert surface.

Additionally, in \cite{AS} the search for further examples was narrowed
through a proof that there are no totally geodesic Seifert surfaces in
two-bridge knot complements.  With the same goal of limiting the existence
of totally geodesic Seifert surfaces in mind, we have the following
theorem of Section \ref{S:topological},

\begin{thm_uniqueness}
Given a semifree totally geodesic Seifert surface $S$ embedded in the
complement of a knot or link $L$, there exists no
other totally geodesic Seifert surface embedded in $S^3 - L$ with the same
boundary slope on each component of $L$.
\end{thm_uniqueness}

Indeed, knowing that all Seifert surfaces in knot complements have the
same boundary slope leads us to the following corrollary:

\begin{cor_uniquenesscor}
Given a semifree totally geodesic orientable Seifert surface $S$ embedded
in the complement of a knot $K$, there exists no other totally geodesic
orientable Seifert surface embedded in $S^3 - K$.
\end{cor_uniquenesscor}

With a similar goal in mind, Section \ref{S:width} defines the
\textbf{width} invariant for surfaces, which is itself motivated by a
definition of width for an $l$-curve.

\begin{defn}
Given a nontrivial, minimal length closed curve $\gamma$ on a maximal
cusp, we call the length of the shortest path which starts and ends on
$\gamma$, but which is not isotopic into $\gamma$, the \textbf{width} with
respect to $\gamma$, sometimes denoted $w_\gamma$.  We sometimes discuss
the width of a knot $K$, denoted $w(K)$, by which we mean the width with
respect to the longitude of $K$.
\end{defn}

\begin{defn}
Let $S$ be a Seifert surface or a non-orientable Seifert surface in the
complement of a hyperbolic knot or link $L$.  Then, by definition,
$\partial S$ is a union of $l$-curves, with exactly one $l$-curve on each
cusp.  We can expand the cusps while forcing the widths of these
$l$-curves on each cusp to remain equal, until there is a cusp 
tangency.  We call the
resulting width the \textbf{balanced width} of the surface $S$. Note that
by definition, the width of a knot must be balanced.
\end{defn}

It turns out that the balanced width of a totally geodesic surface has a
very predictable behavior, leading to the following series of theorems.

\begin{thm_widthlessthantwo}
Consider a hyperbolic knot or link $L$.  If there exists a semifree
totally geodesic Seifert surface $S$, orientable or non-orientable, with
balanced width $w$ in $S^3 - L$, then $w <  2$.
\end{thm_widthlessthantwo}

This bound is actually the best possible, as demonstrated by the free
totally geodesic surfaces in the $(p,p,p)$-pretzel knots (Example
\ref{ex:333}), which have width approaching $2$ from below as $p$
approaches infinity.  On the other hand, the semifree restriction is
indeed necessary, as shown by the totally knotted totally geodesic surface
in Example \ref{ex:totallyknotted} which has width greater than $2$.

\begin{thm_wge1}
Consider a hyperbolic knot or link $L$.  If there exists an embedded
totally geodesic Seifert surface $S$, orientable or non-orientable, with
balanced width $w$ in $S^3 - L$, then $w \ge 1$.
\end{thm_wge1}

Indeed, knowing that $w(S) \ge 1$ for any totally geodesic Seifert surface
allows us to eliminate a very large class of knots from having totally
geodesic Seifert surfaces via the following theorem.

\begin{thm_knotsfan}
Consider a hyperbolic knot $K$ in a reduced, oriented
projection $P$.  Form a sequence of knots $\{K_{i}\}$ by twisting
similarly oriented strands incident to the same region in the projection plane
about each other so as to add an even number of crossings, as
in Figure \ref{fig:tangle}.  If, for some $N > 0$, $n>N$ implies $K_{n}$
is hyperbolic, then $$\lim_{i \to \infty} w(K_{i}) = 0.$$
\end{thm_knotsfan}

\begin{cor_knotsfancor}
With $K_n$ as above, for some positive integer $N$ and all $n > N$, the complement of $K_n$ does not possess a totally geodesic Seifert surface.
\end{cor_knotsfancor}

Note that in the case that the initial projection is a reduced prime
alternating projection that does not correspond to a 2-braid knot, and we twist
to create alternating knots,  all of the
knots will be hyperbolic by results of Menasco. And if the twist makes the
resulting projection nonalternating, it is still true that for enough
twists, the resulting knots will all be hyperbolic.

The final theorem regarding width requires a technical
definition that will be useful throughout the paper.

\begin{defn}
Let $S$ be a semifree surface with boundary in the complement of a link
$L$ and let $D$ be a compressing disk for $\partial(S^3 - N(L) -
N(S))$.  Since $S$ is itself incompressible and boundary incompressible,
$\partial D$ alternates between $n$ arcs in $S$ and $n$ arcs in the cusp
boundaries for some $n > 1$.  Then if $n$ cannot be reduced through
isotopy while preserving the property that $D$ is a compressing disk, we
say that $D$ is an \textbf{essential $n$-gon} in the complement of $S$.
\end{defn}

\begin{thm_we1means3gon}
Let $L$ be a hyperbolic knot or link and $S$ a totally geodesic Seifert
surface, orientable or non-orientable, embedded in $S^3 -
L$ with balanced width $w$.  Then $w=1$ if and only if there is an
essential $3$-gon in the complement of $S$.
\end{thm_we1means3gon}

Finally in Section \ref{S:sixtheorem}, we look at an application of the
examples. The Six Theorem, proven independently by Ian Agol 
\cite{Agol} and Mark Lackenby \cite{Lackenby}, shows that for a finite volume hyperbolic $3$-manifold $N$ with single embedded horocusp $C$, performing Dehn surgery
on a curve $\alpha$ such that the length of $\alpha$ is strictly greater
than six, always yields a hyperbolike manifold. (See Section \ref{S:sixtheorem}
for more details.) Moreover, Agol demonstrated that
this bound is sharp by giving an explicit example of a hyperbolic $3$-manifold
and a curve in its cusp boundary of length exactly six such that Dehn surgery 
on the curve yielded a non-hyperbolike manifold. His example was not a knot complement in the 
3-sphere. In this section, we prove:

\begin{thm_sixtheoremsharp}
The Six Theorem is sharp for knot complements in the 3-sphere, with 
$(p, p, p)$ pretzel knots as examples, for every odd $p \geq 3$.
\end{thm_sixtheoremsharp}

\section{Generating Totally Geodesic Seifert Surfaces}\label{S:generating}

\subsection{Background on 2-Orbifolds and 3-Orbifolds}
An $n$-orbifold is a Hausdorff space $X^n$, along with neighborhoods
locally modelled on $R^n / \Gamma$ where $\Gamma$ is a finite group action.  We
define the singular set of an orbifold to be the set of points in $X^n$
that are locally modelled on $R^n/\Gamma$ where $\Gamma$ is not the identity.
Specifically, the singular set of a 2-orbifold may contain the following:\\

\begin{itemize}
\item Cone points of order n - modelled on $R^2 /Z_n$, where $Z_n$
acts by rotations,\\
\item Corner reflectors of order n - modelled on $R^2/D_n$,
where $D_n$ is the dihedral group of order n, and\\

\item Mirrors - modelled on $R^2/Z_2$, $Z_2$ acts by reflection.\\

\end{itemize}

Likewise, the singular set for a orientable 3-orbifold consists of a
trivalent graph.
The edges of order $n$ are modelled on $R^3/Z_n$, where $Z_n$ acts by
rotations. We label each such edge by $n$ in the orbifold, except when the
edge is modelled on $R^3/Z_2$. The vertices are modelled on the quotient
of $R^3$ by either the dihedral group of order $2n$, the tetrahedral group,
the octahedral group, or the icosahedral group. Thus the edges emanating from
a vertex must be one of the following combinations:  (2, 2, n), (2, 3, 3),
(2, 3, 4), or (2, 3, 5).  For more details on orbifolds, see \cite{Thurston} or
\cite{Dunbar}.

In this paper we will denote 2-orbifolds as $X^2( ; )$, where $X^2$ is the
underlying Hausdorff space, the numbers before the ``;" are cone points
and numbers after the ``;" are the corner reflectors.  In our notation all
points in the
boundary of $X^2$ which are not corner reflectors are mirror points.

\begin{figure}[htbp]
\begin{center}
\includegraphics*[height=1 in]{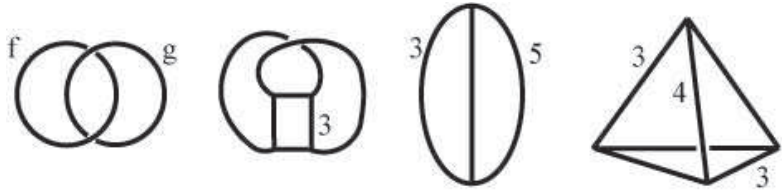}
\end{center}

\caption{\label{fig:3orbifolds}{A few examples of spherical 3-orbifolds.
On the left, $f, g \in Z^+$ }}
\end{figure}

We are specifically interested in spherical 3-orbifolds and rigid
2-orbifolds. A
spherical 3-orbifold is an orbifold with an orbifold covering map from
$S^3$.  The spherical 3-orbifolds are partitioned into a finite number of classes
and the complete list of these classes can be found in \cite{Dunbar}. Figure
\ref{fig:3orbifolds} contains examples of spherical 3-orbifolds.

A 2-orbifold is rigid if it is hyperbolic and its Teichm\"{u}ller space
has dimension
0.  In other words, the orbifold has a unique hyperbolic structure. In a
hyperbolic
2-orbifold, the dimension of the Teichm\"{u}ller space is given by the
function \begin{math}-3\chi(X^2)+2k+l\end{math} \cite{Thurston},  where
$\chi(X^2)$ is
the Euler characteristic of the underlying space, $k$ is the number of
cone points, and $l$ is the number of corner reflectors.  The orbifolds in
Table \ref{2orbi} are the only 2-orbifolds for which the Teichm\"{u}ller
space has
dimension 0.

\begin{table}[h]
\caption{\label{2orbi}{Table of Rigid 2-Orbifolds}}
\begin{tabular}{|l|l|}
\hline
Hyperbolic Rigid 2-Orbifolds&Exceptions (these are not hyperbolic)\\
\hline
\hline
$S^2(n,m,p)$&$S^2(2,2,n), S^2(3,3,3)$\\

\hline
$D^2(n;m)$&$D^2(3;2), D^2(3;3), D^2(4;2)$\\

\hline
$D^2( ;n,m,p)$&$D^2( ;2,2,n), D^2( ;2,3,3), D^2( ;2,3,4),D^2(
;2,3,5),$\\
&$  D^2( ;2,3,6), D^2( ;2,4,4),D^2( ;3,3,3)$\\

\hline
\end{tabular}
\end{table}

Note that a cone point labeled with a positive integer $n$ corresponds to
an elliptic isometry of order $n$.  A cone point labeled with ``$\infty$"
correspond to a parabolic  isometry and can be thought of as a puncture in the interior
of the 2- orbifold. Similarly, a corner reflector point labeled ``$\infty$" is
thought of as a puncture on the boundary of the 2-orbifold.  We define the infinity set of
a hyperbolic 2-orbifold to be the set of infinity cone points and corner
reflectors. For instance, the 2-orbifold $S^2(2, 3, \infty)$ is equivalent to an
open disk (a sphere with a puncture) with cone points of order 2 and 3 and its infinity
set is the boundary of this disk.

\subsection{Embedding of 2-Orbifolds inside 3-Orbifolds}

\begin{thm}Let $J$ be a collection of disjoint arcs and simple closed
curves in a spherical 3-orbifold $N$, such that their complement is a
hyperbolic 3-orbifold $Q$ containing a rigid 2-orbifold $O$ with non-empty infinity set.  If the preimage of $J$ in the covering of $N$ by $S^3$ is a knot or link and if the preimage of $O$
is a Seifert surface $S$ for that knot or link then $S$ is isotopic to a 
totally geodesic Seifert surface.\end{thm}

This appears as Corollary 2.2 in \cite{AS}.

To apply this theorem, we need to consider rigid 2-orbifolds with a
nonempty infinity set. Since we consider a sphere with one ``$\infty$" cone point to
be the same as an open disk, we think of a puncture from an ``$\infty$" cone
point in a 2-orbifold as removing a closed disk from the interior of the 2-orbifold.  We will often represent  ``$\infty$" cone points as ``$\infty$"
closed loops.  Likewise, an ``$\infty$" corner reflector can be thought of
as removing a closed disk, centered at the corner reflector, from the
2-orbifold.  We will often represent an ``$\infty$" corner reflector as an ``$\infty$" arc
in the 2-orbifold.

\begin{figure}[htbp]
\begin{center}
\includegraphics*[height=2 in]{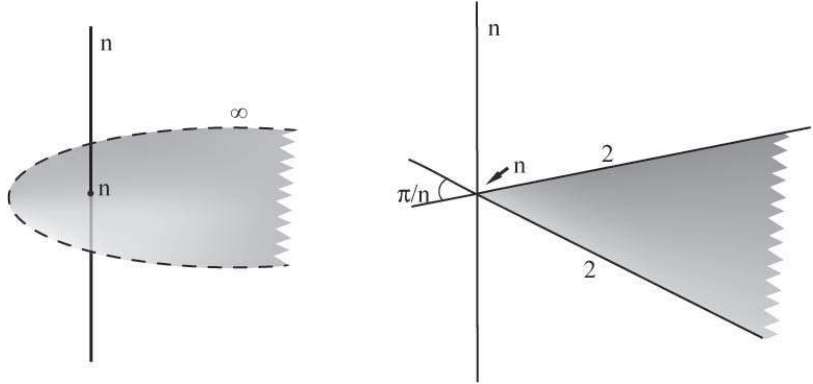}
\end{center}
\caption{\label{fig:2orbs}Left: cone point of order $n$. Right: Corner
reflector of order $n$.}
\end{figure}

Because we want the preimage of the infinity set in $S^3$ to be a link,
any ``$\infty$" arc in a 2-orbifold must end on a 2-axis of the
3-orbifold. When an
orbifold has two ``$\infty$"  corner reflectors with a path connecting
them that consists of only mirror points, then their corresponding ``$\infty$" arcs
must end on a common 2-axis. An orbifold with three infinity corner reflectors,
will have three arcs where each pair of arcs end on a common 2-axis. Since $D(n;
\infty)$ has only one corner reflector, its ``$\infty$" arc must start and end on
the same axis. Cone points that are not in the infinity set, say of degree $n$, are
realized in the 2-orbifold by an intersection with an axis of order $n$ in the
3-orbifold.  Corner reflectors of order $n$ are realized in the 2-orbifold as the
intersection point of two 2-axes in the 3-orbifold. Examples of these can be seen in
Figure \ref{fig:2orbs}. The order of the corner reflector corresponds to
the angle $\pi/n$ between the 2-axes.

Finally, we must ensure that the resulting link is hyperbolic. The
creation of essential tori, annuli and spheres must be avoided.

Immersed totally geodesic surfaces can also be generated from a
similar process using immersed rigid 2-orbifolds.  Examples of immersed
surfaces appear in the following examples section.

\subsection{Examples}\label{SS:examples}

Now that we have the background, we can look at a few interesting examples
of totally geodesic surfaces generated with this method.

\begin{figure}[htbp]
\begin{center}
\includegraphics*[height=2 in]{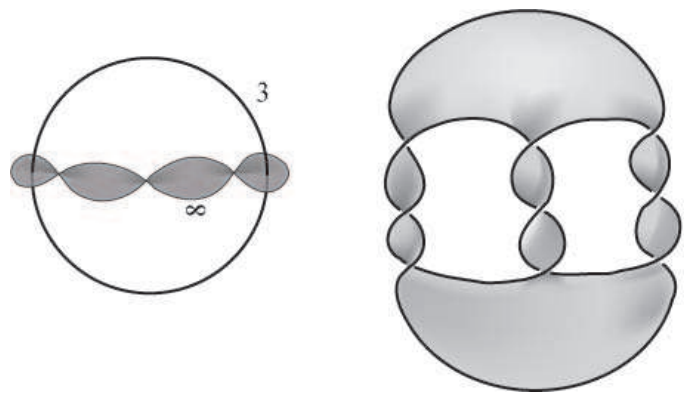}
\end{center}

\caption{\label{fig:333}The $(3, 3, 3)$ pretzel knot, and a generating
orbifold.}
\end{figure}
\begin{figure}[htbp]
\begin{center}
\includegraphics*[height=2 in]{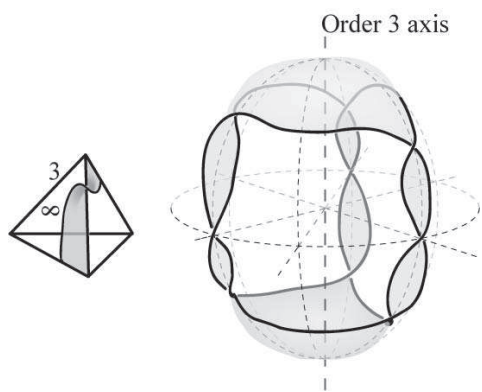}
\end{center}

\caption{\label{fig:333a}Another view of the $(3, 3, 3)$ pretzel knot with
another
generating orbifold.}
\end{figure}

\begin{example}\label{ex:333} The $(3, 3, 3)$ pretzel knot.
\end{example}
A $(p, p, p)$ pretzel knot is a knot with 3 arms, each of which contains p
crossings. The $(3, 3, 3)$ pretzel knot is shown in Figure \ref{fig:333} with a
generating orbifold. The grey surface area in the figure is a totally geodesic
Seifert surface in the knot complement.  Note that there
are multiple ways of embedding a rigid 2-orbifold in a spherical
3-orbifold to lift to the $(3, 3, 3)$ pretzel knot: in Figure \ref{fig:333} the rigid
2-orbifold is a $S^2(\infty, 3, 3)$, while in Figure \ref{fig:333a} it
is a $D^2( ; \infty, 2, 3)$. In Figure \ref{fig:333a}, the knot is drawn
with symmetry axes corresponding to the axes of
the generating 3-orbifold.

Similar surfaces can be made in any $(p, p,\dots,p)$ pretzel knot complement.

\begin{figure}[htbp]
\begin{center}
\includegraphics*[height=2 in]{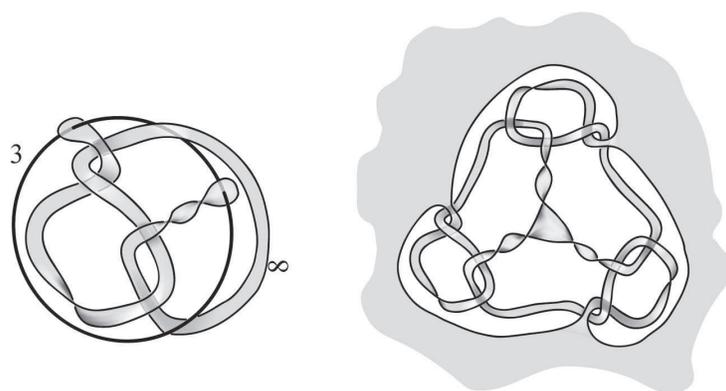}
\end{center}

\caption{\label{fig:totally}A totally knotted surface.}
\end{figure}

\begin{example} A totally knotted surface.\label{ex:totallyknotted}
\end{example}
Another way to make more complicated knots is to knot up the
$S^2(\infty, 3,
3)$ orbifold, as in Figure
\ref{fig:totally}.
The result is a totally knotted totally geodesic surface. This surface, as
in the
previous example, is orientable. It is of note that this knot has width
greater than
2; thus 2 as an upper bound on width for semifree surfaces does not hold
in the
totally knotted case.

\begin{figure}[htbp]
\begin{center}
\includegraphics*[height=2 in]{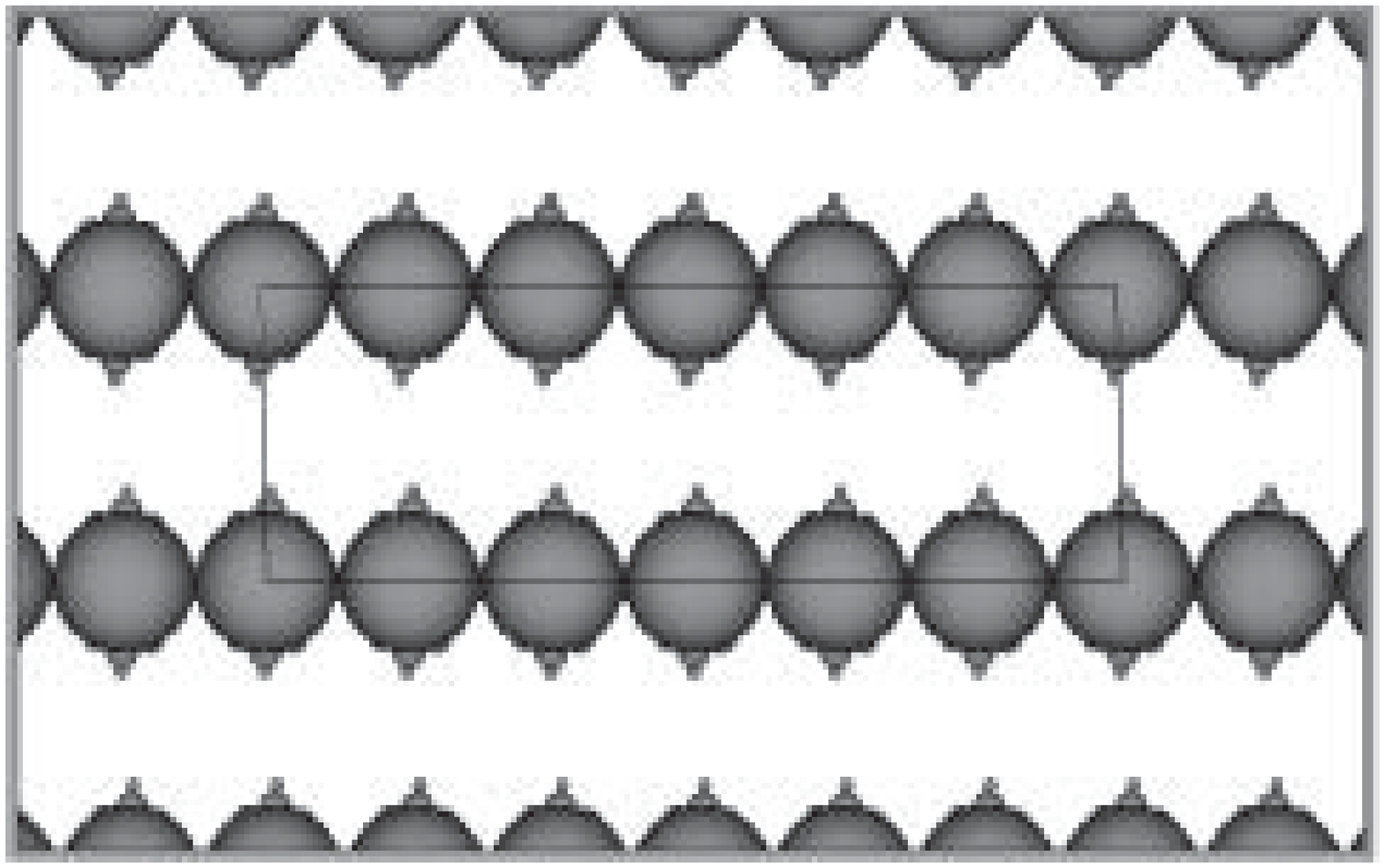}
\end{center}

\caption{\label{fig:totallyhoroball}The horoball diagram for the totally
knotted
surface.}
\end{figure}

\begin{figure}[htbp]
\begin{center}
\includegraphics*[height=2 in]{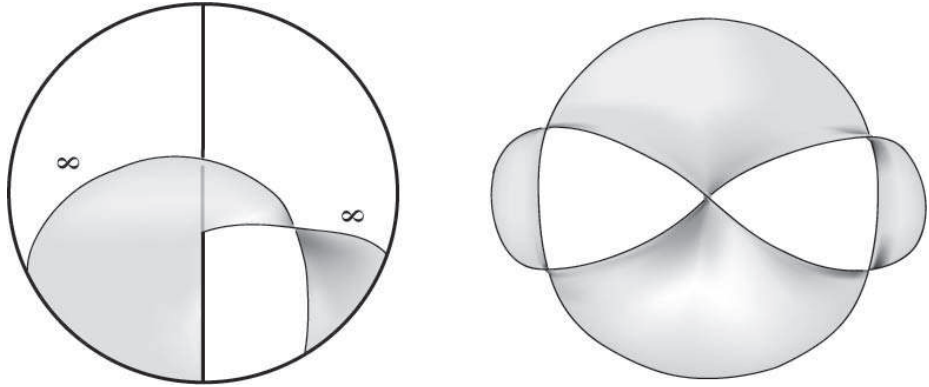}
\end{center}

\caption{\label{fig:whitehead}The Whitehead link.}
\end{figure}

\begin{example} The Whitehead link.
\end{example}
An example of a nonorientable totally geodesic checkerboard surface in a
link
complement is found in the Whitehead link, shown in Figure
\ref{fig:whitehead}. It is
not known whether nonorientable totally geodesic Seifert surfaces exist in
knot
complements.

By twisting up the two infinity arcs, or by twisting one arc around the
2-axis it
ends on, we can create a family of links containing totally geodesic
surfaces.
These are (2p, 2q + 1, 2p) pretzel links. These orbifolds will always lift
to links,
and always produce nonorientable totally geodesic surfaces.

\begin{figure}[htbp]
\begin{center}
\includegraphics*[height=2 in]{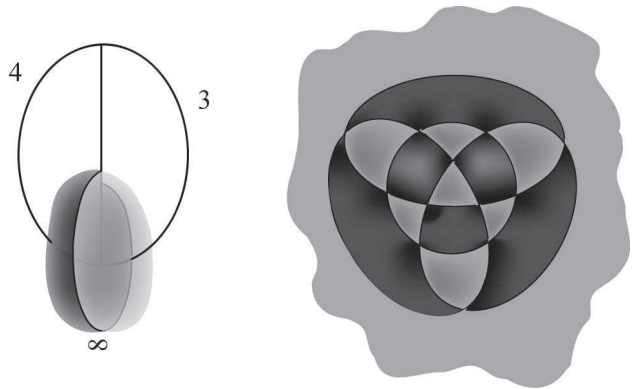}
\end{center}

\caption{\label{fig:34link}A link with two embedded totally geodesic
checkerboard surfaces.}
\end{figure}

\begin{example}\label{ex:34link} Multiple totally geodesic surfaces.
\end{example}
The link in Figure \ref{fig:34link} is of particular interest because it
has two totally
geodesic checkerboard surfaces. The link can be realized as the lift of
two rigid
orbifolds at once--a $D^2(3; \infty)$ and a $D^2(4; \infty)$. Two other
links are
known to have two totally geodesic checkerboard surfaces; they are
generated
by a similar configuration: by orbifolds sitting in a spherical 3-orbifold
with axes
labelled (2, 3, 3) (this gives the Borromean rings) and (2, 3, 5).

\begin{figure}[htbp]
\begin{center}
\includegraphics*[height=2 in]{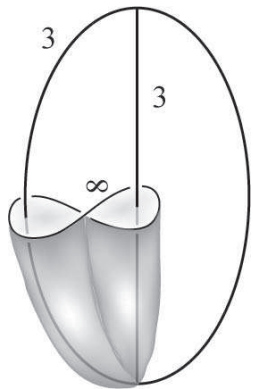}
\end{center}

\caption{\label{fig:33orb}Another orbifold that lifts to this link.}
\end{figure}

There is also an immersed totally geodesic surface in this link
complement. The
surface that comes from $D^2(3; \infty)$ can also be lifted from the
$S^2(\infty,
3, 3)$ shown in Figure \ref{fig:33orb}. The surface in grey is a self-
intersecting $S^2(2, \infty, \infty)$. It lifts to a set of four
thrice-punctured disks,
all intersecting each other, each with with a different link component as
boundary.

\begin{figure}[htbp]
\begin{center}
\includegraphics*[height=2 in]{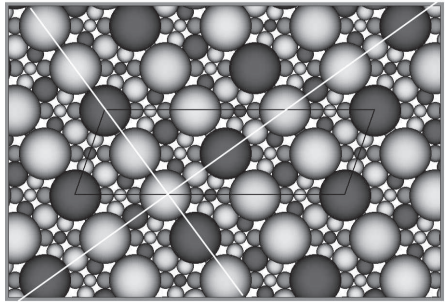}
\end{center}

\caption{\label{fig:34horoballs}Horoball diagram.}
\end{figure}

All three of these totally geodesic surfaces can be seen in the horoball
diagram in Figure \ref{fig:34horoballs}. The immersed surface runs along the
longitude and meridian, while the embedded surfaces run along lines drawn in white.

\section{Using Topological Means to Disprove the Existence of Totally
Geodesic Surfaces} \label{S:topological}

There are some properties, both in knot and link
complements and in the surfaces themselves, that allow us to eliminate
surfaces from contention as totally geodesic candidates, and 
eventually to prove a fact about the uniqueness of totally goedesic surfaces in a given knot or link complement. In \cite{AS}, the authors use the topological properties of totally geodesic surfaces to show that no such
surface can have a bigon in its complement.

\begin{thm}\label{bigons}
Any surface $S$ with an essential bigon in its complement cannot be
totally geodesic.
\end{thm}

\begin{proof} Cut the manifold open along the totally geodesic
surface and double it.  The essential bigon doubles to an essential annulus, contradicting the hyperbolicity of
the doubled manifold.
\end{proof}

In the case of a checkerboard surface for an alternating knot or link, these
bigons can occur in any reduced alternating projection.

\begin{thm}\label{n-gons}
An n-gon region $R$ in the projection plane of a reduced alternating
diagram in the complement of a totally geodesic checkerboard surface $S$
must correspond to an essential n-gon.
\end{thm}

\begin{proof}
Each of the crossing arcs corresponding to this projection are essential
arcs in the surface $S$. Hence, since $S$ is totally geodesic, they 
must lift to geodesics, each connecting two distinct horoballs in {\Hth}. Since 
$R$ must lift to a collection of disks in {\Hth}, it lifts to a disk that is bounded 
by an alternating sequence of $n$ horoballs and geodesics connecting the horoballs. Hence $R$ is essential.
\end{proof}

Now, we may observe that in a given projection $P$ of a link $L$,
checkerboard surfaces can only be totally geodesic if they contain no bigons in their
complement in $P$. Since checkerboard surfaces are complementary in $P$,
any $P$ for which both checkerboard surfaces contain bigons cannot have a
totally geodesic checkerboard surface. Note that Example \ref{ex:34link}
yields an example of an alternating link with no bigons in a 
projection and both checkerboard surfaces totally geodesic.

\subsection{Uniqueness of Totally Geodesic Seifert Surfaces for Knots and
Links}

In many cases, if there exists a totally geodesic Seifert surface, it is
unique.

\begin{thm}\label{uniqueness}
Given a semifree totally geodesic Seifert surface $S$ embedded in the
complement of a knot or link $L$, there exists no other totally geodesic
Seifert surface embedded in that complement with the same boundary slope
on each component of $L$.
\end{thm}

\begin{proof}
Since $S$ is semifree, there exists a compressing disk $D$. Denote the
pre-image of $S$ in {\Hth} as $p^{-1}(S)$ and the pre-image of D as
$p^{-1}(D)$. For
a particular copy of $D$ in $p^{-1}(D)$, called $\widetilde{D}$, the
boundary of $\widetilde{D}$ is a curve which lies alternatingly on a cyclic sequence of
geodesic planes in $p^{-1}(S)$ and the series of horoballs that occur at their
points of tangency. Let that cyclic sequence of geodesic planes be denoted
$\widetilde{S_1}$, \dots , $\widetilde{S_n}$, where $n \ge 3$. Denote 
the horoball that occurs at the
point of tangency between $\widetilde{S_i}$ and $\widetilde{S_{i+1}}$ as $H_i$, and
call that point $P_i$.  We will consider $H_1$ to be the
horoball at infinity. (See Figure ~\ref{fig:BigPicture}.)

\begin{figure}[htbp]
\begin{center}

\includegraphics*[height=1.5in]{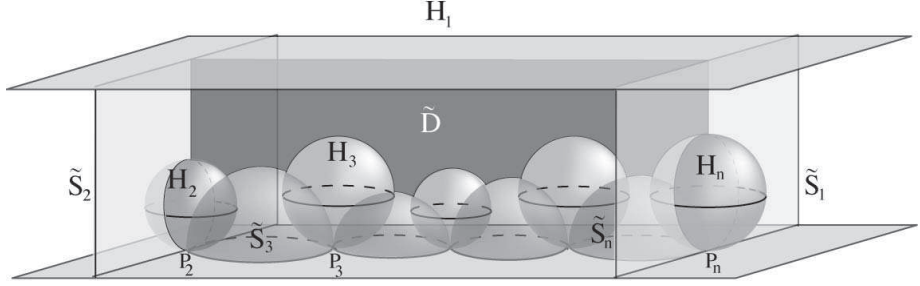}
\end{center}

\caption{\label{fig:BigPicture} The compressing disk $\widetilde{D}$ 
borders on a series of geodesic planes and horoballs.}
\end{figure}

Now, assume there exists another totally geodesic Seifert surface $T$
embedded
in $S^3 - L$ which has the same boundary slope as $S$ for each component
of $L$, but is distinct from $S$. A meridian of $L$ thus intersects the
boundaries of $T$ and $S$ exactly once each. As we continue to trace the meridian
multiple times, we alternate between intersections with the boundary of $T$ and
intersections with the boundary of $S$. Thus, between the boundary of any
two lifts of S in $p^{-1}(S)$, $\widetilde{S_{a_1}}$ and
$\widetilde{S_{a_2}}$,
there must exist a geodesic plane $\widetilde{T'}$ in the pre-image $p^{-
1}(T)$. For any two lifts of $S$ tangent at the base of some $H_i$, then,
there must be a lift of $T$ which separates them.

Recall that $H_1$ is the horoball at infinity. Then, $\widetilde{S_1}$ and
$\widetilde{S_2}$ take the form of vertical geodesic planes, with
$\widetilde{S_3}$ a hemispherical plane tangent to $\widetilde{S_2}$ at
$P_2$.  Recall that there must be a lift of $T$ in $p^{-1}(T)$ which separates
each pair $\widetilde{S_i}$ and $\widetilde{S_{i+1}}$, and thus is tangent to
them both at $P_i$. Any vertical geodesic plane $\widetilde{T_1}$ tangent to
$P_2$ will be identical to $\widetilde{S_2}$, contradicting our assumption that
$S$ is distinct from $T$; by the same logic, $\widetilde{T_1}$ cannot be the
hemispherical geodesic plane tangent to both $P_2$ and $P_3$, as it would
then be identical to $\widetilde{S_3}$. Lastly, $\widetilde{T_1}$ cannot be a
hemispherical plane tangent to $\widetilde{S_3}$ at $P_2$ with smaller
diameter than $\widetilde{S_3}$, or it would exist entirely inside
$\widetilde{S_3}$,
and not between it and $\widetilde{S_2}$. Thus, $\widetilde{T_1}$ must be
a hemispherical geodesic plane tangent to $\widetilde{S_3}$ at $P_2$ with
diameter greater than that of $S_3$. Thus, $\widetilde{T_1}$ contains both
$P_2$ and $P_3$. But, since there must exist some lift of $T$ in between
every pair of tangent copies of $S$, there must be some $\widetilde{T_2}$
tangent to $P_3$, which necessarily intersects $\widetilde{T_1}$. Hence, $T$ is not an
embedded surface in $S^3 - L$. Thus, we again have a contradiction, and
our assumption of the existence of $T$ must be false.
\end{proof}

In the case that $L$ is a knot, and $S$ is orientable, we have the
following corollary:

\begin{cor}\label{uniquenesscor}
Given a semifree totally geodesic orientable Seifert surface $S$ embedded
inthe complement of a knot $K$, there exists no other totally geodesic
orientable Seifert surface embedded in $S^3 - K$.
\end{cor}

\begin{proof}
Since $S$ is an orientable Seifert surface in a knot complement, it has
boundary slope parallel to the longitude. Any other such surface must also
have the same boundary slope. By the previous theorem, there can be no such
surface distinct from S.
\end{proof}

\section{The Width Invariant for Totally Geodesic Surfaces}\label{S:width}

In this section, we consider how width can impact the possible existence
of totally geodesic Seifert surfaces.

\subsection{Bounds on Width}

\begin{thm}\label{widthlessthantwo}
Consider a hyperbolic knot or link $L$.  If there exists a semifree
totally geodesic Seifert surface $S$, orientable or non-orientable, with
balanced width $w$ in $S^3 - L$, then $w <  2$. This
upper bound is best possible.
\end{thm}

\begin{proof}
Since $S$ is semifree, there is a compressing disk $D$ in $S^3
- (N(K) \cup N(S))$.  The boundary of $D$ consists of arcs which
alternate between lying in the boundary of the cusp set $\{C_i\}$ and
lying in the
surface.  Indeed, the set $\partial D \cap (\cup \{C_i\})$ is a collection
of arcs each of which travels on some element of $\{C_i\}$ nontrivially
from $S$ back to $S$.  Thus, the length of each arc in this set is
greater than or equal to the balanced width $w$ of the cusp set $\{C_i\}$.

The disk $D$ lifts to a collection of closed disks in {\Hth} - let
$\widetilde{D}$ be one such copy.  Then $\partial \widetilde{D}$
alternates between travelling along the boundaries oof horoballs covering the
cusp set and geodesic planes covering the surface. Choose one such
horoball to be centered at $\{\infty\}$, similar to Figure 
\ref{fig:BigPicture}.

Since there are only a finite number of horoballs in this chain, there
must be some horoball $A$ with Euclidean height less than or equal to the Euclidean height of every other horoball in the chain, excluding the horoball at infinity.  But
$\partial \widetilde{D} \cap A = \gamma$ is a curve which starts and ends
at points on A that are tangencies with horoballs no smaller than A. 
Hence $\gamma$ starts and ends
at or above the equator of $A$.  Using the triangle inequality and 
the fact that the distance from the top of a horoball to the equator is always exactly
$1$, we see that $|\gamma| \le 2$.  But $w \le |\gamma|$ and so we see
that $w \le 2$.

Now consider the case where $w = 2$.  Then $|\gamma|=2$ and since no other
horoball can be strictly smaller than $A$ we see that the horoballs on
either side of $A$ in the sequence are actually both the same height as
$A$ and tangent to $A$.  In fact, we are forced to have a sequence of
equal height, tangent horoballs.  But now consider a horoball $B$ adjacent
in the sequence to the horoball at $\{\infty\}$.  $\partial \widetilde{D}
\cap B$ is a curve which starts at the equator of $B$ and ends at the top
of $B$, forcing $w$ to be less than or equal to $1$, contradicting the
assumption that $w=2$.

The $(p,p,p)$ pretzel knots of Example \ref{ex:333} yield a sequence 
of hyperbolic knots with
width approaching 2 from below and with
free(and hence semifree) totally geodesic Seifert surfaces, demonstrating
that the upper bound of 2 is best possible.
\end{proof}

\begin{thm}\label{wge1}
Consider a hyperbolic knot or link $L$.  If there exists an embedded
totally geodesic Seifert surface $S$, orientable or non-orientable, with
balanced width $w$ in $S^3 - L$, then $w \ge 1$.
\end{thm}

\begin{proof}
Assume $w < 1$. The totally geodesic surface $S$ has
boundary a union of $l$-curves. We maximize the cusps
while forcing the widths with respect to these
$l$-curves to be equal. Hence not every cusp will
necessarily have a point of tangency. Let $C$ be a
cusp with a point of tangency with itself or another
cusp and let $\widetilde{C}$ be a horoball covering
$C$ centered at $\{ \infty \}$, and normalized to have boundary a
horizontal plane of Euclidean height 1. There
is a horoball $A$ tangent to
$\widetilde{C}$.

  The surface $S$ lifts to a set of geodesic planes
containing two vertical planes that are a distance $w$
apart. The width curve with respect to the resulting
$l$-curves on the horosphere at infinity has a well
defined direction. If we travel along a great circle
in this direction from the top of $A$ a distance at
most $w$, we will have reached a hemisphere
$\widetilde{S}$ contained in the lift of $S$. Since $w
< 1$ and the hyperbolic distance from the top of $A$
to the equator is $1$, $\widetilde{S}$ must intersect
$A$ above the equator, hence $\widetilde{S}$ has
radius that is greater than $\frac{1}{2}$. The surface
is embedded, thus $\widetilde{S}$ is contained between
two vertical planes contained in the lift of $S$. The
distance between the two vertical planes is equal to
$w$, both in the Euclidean and hyperbolic length since
the maximal cusp is normalized to height $1$.  But
since the radius of $\widetilde{S}$ is greater than
$\frac{1}{2}$, this implies that $w > 1$, which
contradicts the assumption.

\end{proof}

\begin{thm}\label{we1means3gon}
Let $L$ be a hyperbolic knot or link and $S$ a totally geodesic Seifert
surface, orientable or non-orientable, embedded in $S^3 -
L$ with balanced width $w$.  Then $w=1$ if and only if there is an
essential $3$-gon in the complement of $S$, and this can occur only if $S$
is nonorientable.
\end{thm}

\begin{proof}
First, assume there is an essential $3$-gon $D$ in the
complement. Note that this implies $S$ is nonorientable, since otherwise,
let $S_{+}$ and $S_{-}$ be the two copies of $S$ on the boundary of the
regular neighborhod of $S$. Arcs in $\partial D \cap \partial N(S)$
must alternate between lying in $S_{+}$ and $S_{-}$. Hence there must be
an even number of them.

The essential $3$-gon $D$ lifts to a disk $\widetilde{D}$ in {\Hth}
bounded by two vertical planes $V_1$ and $V_2$ and a hemisphere $\widetilde{S}$
covering the totally geodesic surface and three horospheres covering the cusp boundary, one of
which is centered at $\infty$ and is denoted $\widetilde{C}$ and the other two of
which are denoted $A$ and $B$. (See Figure ~\ref{fig:wisone}.)

\begin{figure}[htbp]
\begin{center}
\includegraphics*[height=2.3in]{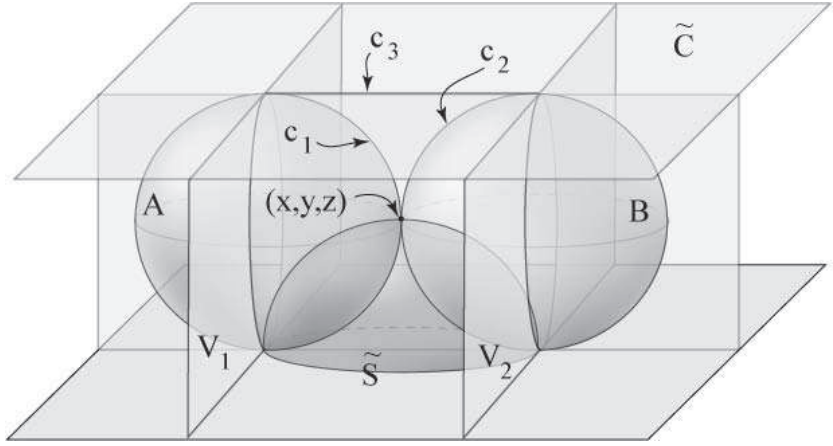}
\end{center}
\caption{\label{fig:wisone} The case when an essential 3-gon is present.}
\end{figure}

  The hemisphere $\widetilde{S}$ meets
each vertical plane only at the center of $A$ and $B$
since the totally geodesic surface is embedded. The
$3$-gon $D$ has arcs in its boundary that lie in the
surface. These boundary curves cannot be isotoped to
the surface since the number of boundary curves of $D$ would
then not be minimal and so $D$
would not be essential. Hence for the lifts, $c_i$, of
the boundary curves, $ |c_i| \ge w$  for all $i$. Let
the origin of the upper-half space model of {\Hth} be
taken as the center of $A$ on the boundary of
the $xy$ plane such that $A$ and $\widetilde{S}$ are
centered on the $y$ axis. Simple calculations show
that for the height, $z$, of the point of
intersection, $(x,y,z)$, of $A$ with $\widetilde{S}$
and the $yz$ plane, $z = \frac{r}{a}x$, where $r$ is
the radius of $\widetilde{S}$ and $a$ is the radius of
$A$. Since $\widetilde{C}$ is at Euclidean height $1$,
$2r \ge w$. Therefore $\frac{1}{2} \le r \le 1$, and
$0 < a \le \frac{1}{2}$. Hence $z \ge
a$. Since $w$ is also realized as a segment of a great
circle running from the top of $A$ to $(x,y,z)$, it
follows from the lower bound on $z$ that $w \le 1$,
since the segment of a great circle on $A$ to the
equator of $A$ is $1$ in hyperbolic distance and
$(x,y,z)$ is above or at the equator of $A$. Theorem
$\ref{wge1}$ then implies that $w=1$.

\medskip{}

Now, assume there is an embedded totally geodesic
surface $S$ and $w=1$. Let $C$ be a cusp with a point of
tangency with itself or another cusp and let
$\widetilde{C}$ be a horoball covering $C$. Center
$\widetilde{C}$ at $\{ \infty \}$ and normalize it to height 1. There is a
  horoball $A$ tangent to $\widetilde{C}$. The surface $S$
lifts to a set of geodesic planes containing two
vertical planes that are a distance $w$ apart. The
width curve with respect to the resulting $l$-curves
on the horoball at infinity has a well defined
direction.

If the cusp does not touch itself in $S$, we may travel along a great 
circle on $A$ in
the well defined direction a distance at most $w$ and
we will have come to a hemisphere. Since we assume a
vertical plane does not intersect the top of $A$ there
must be a hemisphere intersecting $A$ above its
equator, as the hyperbolic distance from the top of
$A$ to the equator is $1$. This implies that the
radius of the hemisphere is greater than
$\frac{1}{2}$. But since $S$ is embedded, the
hemisphere is also contained between two vertical
planes a Euclidean and hyperbolic distance of $w = 1$
apart, which is impossible.

Thus the cusp must touch
itself in the surface $S$. Hence there is a vertical
plane, $V_1$, containing a boundary curve of
$\widetilde{C}$, centered on $A$. If we travel along a
great circle on $A$ in a well defined direction with
respect to the resulting $l$-curve on the horoball at
$\{ \infty \}$ from the top of $A$ a hyperbolic
distance of $1$ we will have reached a geodesic plane,
since the width, $w$, is equal to $1$. Thus there is a
hemisphere, $\widetilde{S}$, intersecting $A$ at
Euclidean height $\frac{1}{2}$. If we travel in the
same direction as the great circle along a straight
line in $\widetilde{C}$ a distance of $1$, we will
again have reached a geodesic plane, hence there is a
second vertical plane, $V_2$, a distance of $1$, both
in Euclidean and hyperbolic distance, from $V_1$. The
surface is embedded, therefore $\widetilde{S}$ does
not intersect $V_1$ or $V_2$ except perhaps at a
point. It follows that the radius of $\widetilde{S}$
is $\frac{1}{2}$, and $V_1$, $V_2$, and
$\widetilde{S}$ then bound a disk that is an essential
$3$-gon in {\Hth} with $V_1$ and $V_2$ intersecting
$\widetilde{S}$ at distinct points on the boundary of
{\Hth} and meeting each other at $\{ \infty \}$. This
disk projects to an essential $3$-gon in the manifold.
\end{proof}

\subsection{An Application of Width}

\begin{defn}

A {\bf minimal $l$-curve} for the maximal cusp $C$ of a hyperbolic knot
$K$ is its $l$-curve of shortest length.

\end{defn}

Under the above definition, a minimal $l$-curve always exists, although
there could potentially be two $l$-curves of shortest length.  We will use
minimal $l$-curves to show that the width for certain knots is small.  (Recall that
when we speak of the width of a knot, without regards to a particular $l$-curve,
we mean width with respect to the longitude.)

\begin{figure}[hbtp]

\begin{center}
\includegraphics*[height=1.5in]{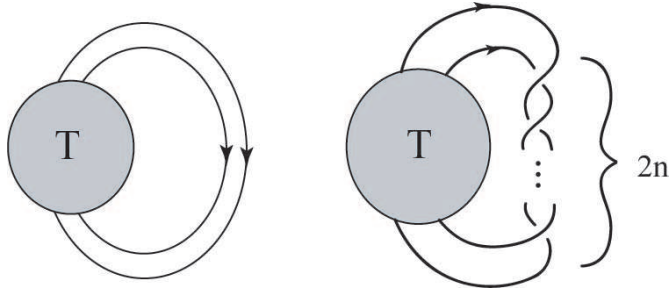}

\caption{\label{fig:tangle}  The second knot is obtained from the first
by adding an even number of crossings to two similarly oriented strands.}
\end{center}
\end{figure}

Consider a knot $K$ in an oriented projection $P$.  Given two similarly
oriented strands of the knot that are both on the boundary of the same
region within the projection, we
can form a sequence of knots $\{K_{p}\}$ by twisting the
strands about each other so as to add an even number of crossings, as in Figure
\ref{fig:tangle}.

\begin{thm}\label{knotsfan}
If, for some $N > 0$, $n>N$ implies $K_{n}$ is
hyperbolic, then $$\lim_{p \to \infty} w(K_{p}) = 0.$$
\end{thm}

\begin{proof}

We first describe a method for obtaining the knots $K_{p}$ described
above. Begin with the knot $K$ and create a link $L$ by drilling out a curve
$\gamma$ which bounds a disk $D$ that is punctured by the two strands of $K$ we
wish to twist, so that the strands have the same orientation as they 
pass through $D$. Then, upon performing $(1, p)$ Dehn filling on $\gamma$, we obtain the
knot $K_{p}$.  Working directly with the link $L$ will help us to show that the
width for these knots becomes small.

From the link $L$, before the Dehn filling, we can form a series of links,
the complements of which are all homeomorphic.  Let $L_{0}$ denote $L$ and let
$L_{p}$ denote the result of adding $2p$ twists to $L_{0}$ in the above
manner.  The complements are indeed homeomorphic, because they are formed by
cutting the complement of $L$ open along the disk $D$, twisting one copy
of $D$ $p$ times, and gluing together again in the original manner.

Let $\alpha$ denote $K$'s longitude
and $\beta$ denote $K$'s minimal $l$-curve.  Let $\alpha_{p}$ and
$\beta_{p}$ denote respectively the images of $\alpha$ and $\beta$ under the
homeomorphisms $h_{p}$ from $L_{0}$ to $L_{p}$.  Mostow's Rigidity Theorem
guarantees that the hyperbolic structures of the link complements are the
same. For this reason, the image of $L_{0}$'s minimal $l$-curve under $h_{p}$
will also be the minimal $l$-curve for $L_{p}$.  However, this does not hold
true for the longitude.  Let $\eta_{p}$ denote the longitude of $L_{p}$.

Assume the intersection number of $\alpha$ with $\beta$ is $x$.
Computations show that the linking number of $\alpha_{p}$ with the core curve of $K$'s
image is $\pm 4p$, where the sign depends on the direction in which the twisting
is done.  (This computation does not work if the strands are oppositely
oriented, in which case the linking number is zero.)  The above core curve has by
definition linking number zero with $\eta_{p}$.  Thus, the intersection number of
$\beta_{p}$ with $\eta_{p}$ is $x \pm 4p$, because homeomorphisms preserve
intersection number.  So, as $p$ approaches infinity, the linking number
of the curve $\beta_{p}$ (which remains constant on the cusp) with $\eta_{p}$
approaches infinity.  Therefore, $| \eta_{p} | \to \infty$ as $p \to
\infty$.  Because the cusp area $A$ must remain constant under the homeomorphisms, and
because $A$ is equal to the product of the longitude length and the 
width, we must have $w_{p} \to 0$, where this width is with respect to $\eta_{p}$.

Performing $(1, 0)$ Dehn filling on the image of $\gamma$ under the map
$h_{p}$, which corresponds to performing $(1, p)$ surgery on the original
$\gamma$, gives a {\it knot} complement which, in general, could have a
quite different hyperbolic structure from \Sth $- L_{p}$.  However, if we take
$p$ large enough, then these structures get arbitrarily close.  In particular, we
can choose $p$ so that the width after the filling also gets arbitrarily close to
zero, as required.

\end{proof}

The above construction of knots $K_{p}$ relies on the fact that all of
them past a certain point are hyperbolic.  One instance where this occurs is
when the knot $K$ is a reduced prime alternating knot that is not a two
braid knot. Then, twisting similarly oriented strands in any reduced alternating
projection so that the resulting knots are alternating will yield a
sequence of  knots, all of which are hyperbolic, by results in 
\cite{Menasco}.  Hence, the knots will have width
approaching zero. Moreover, in the same case, if we twist in the direction
that yields nonalternating knots, the resulting knots will still be
hyperbolic for large enough twists, since they will be limiting toward an
augmented alternating link, which was shown to be hyperbolic in
\cite{Adams}.

Note that if the orientations of the two strands do not
match, the resulting sequence of knots need not have width approaching 0,
as occurs for the sequence of twist knots.

Our main interest in Theorem \ref{knotsfan} lies in the following
corollary:

\begin{cor}\label{knotsfancor}
Consider, as above, a sequence of knots $\{K_{p}\}$ obtained by twisting similarly oriented strands
about each other in a projection of a knot $K$, as in Figure \ref{fig:tangle}.  If, eventually, all knots past a
certain point in the sequence are hyperbolic, then we can find $N>0$ so that
$n>N$ implies the knot $K_{n}$ cannot possess any totally geodesic Seifert
surfaces.
\end{cor}

\begin{proof}
This follows immediately from Theorem \ref{wge1} and Theorem
\ref{knotsfan}, since we can make width arbitrarily small.
\end{proof}

\section{Application: The Six Theorem is Sharp for Knot
Complements}\label{S:sixtheorem}

A manifold is said to be {\it hyperbolike} if it is irreducible with
infinite word-hyperbolic fundamental group.  Under this definition, 
hyperbolic and hyper{\it bolike} manifolds are very similar: for instance, neither can possess an
essential torus.  In fact, a proof of Thurston's geometrization conjecture would
imply that hyperbolic and hyperbolike manifolds are exactly the same.

The Six Theorem, proven independently by Ian Agol \cite{Agol} and Mark
Lackenby \cite{Lackenby}, showed that, for a finite volume hyperbolic $3$-
manifold $N$ with single embedded horocusp $C$, performing Dehn surgery
on a curve $\alpha$ in the cusp boundary such that the length of 
$\alpha$ is strictly greater than six
always yields a hyperbolike manifold.  Moreover, Agol demonstrated that
this bound is sharp by giving an explicit example of a hyperbolic $3$-manifold
and a curve of length exactly six such that Dehn surgery on it yielded a non-
hyperbolike manifold.  In this section, we demonstrate that, furthermore,
the bound is sharp for hyperbolic knot complements, a case that was not
covered by Agol's example.  Our first task is to find a candidate curve on which we
can perform the surgery.

\begin{lem}\label{lowerbound}
For all $(p, p, p)$ pretzel knots, with odd $p \geq 3$, the longitude
length is greater than or equal to six.
\end{lem}

\begin{proof}

Figure \ref{fig:333a} from an earlier section shows symmetries of a $(3, 3, 3)$
pretzel knot.  The $(p, p, p)$ case  with odd $p \geq 3$ is completely
analogous. The vertical axis represents a rotational symmetry of order three, and the
circular axis running horizontally along the equator is a rotational
symmetry of order two.   (There are other symmetries, but these are the
two that will concern us.)  These correspond to isometries of \Hth.
Because these symmetries preserve the totally geodesic surface which has
boundary along a longitude, they must send longitude to longitude.  Since
neither of the symmetry axes touches the knot, they must both correspond
to parabolic isometries.  Combining this information, we have that the
parallelogram corresponding to the fundamental domain of the cusp should
contain symmetries realized as longitudinal translations of order two and
order three.

As usual, consider the horoball at infinity to be normalized so that its
height is one and consider any full-sized horoball tangent to it.  Because the
symmetries preserve the horoball diagram, they force a minimum of six full-sized
horoballs lying along a longitude.  As mentioned in the preceding paragraph, these
symmetries all occur within one fundamental domain of the cusp, and so the
longitude length must have room for all six full-sized balls.  Because
they all have diameter one, this forces the longitude length to be greater than or
equal to six, as desired.
\end{proof}

\begin{figure}[hbtp]

\begin{center}
\includegraphics*[height=1.5in]{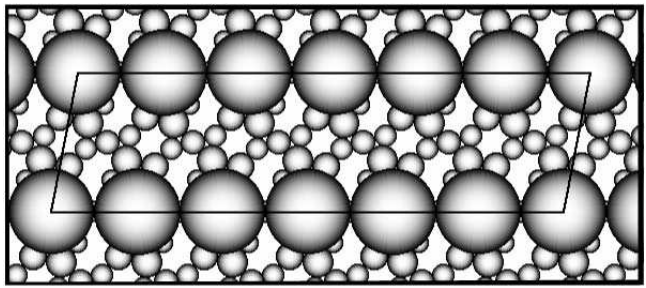}
\caption{\label{fig:horoball}  The horoballs do indeed satisfy an order
six translational symmetry along a longitude.}
\end{center}
\end{figure}

We can see this phenomenon explicitly in Figure \ref{fig:horoball},
provided by SnapPea (see \cite{Weeks}).

One more  lemma is needed before the main result, which will follow
immediately.

\begin{lem}\label{upperbound}
Performing Dehn surgery on the longitude of the $(p, p, p)$ pretzel knot,
for odd $p \geq 3$, yields a non-hyperbolike manifold.
\end{lem}

\begin{proof}

Consider the totally geodesic surface in the $(p, p, p)$ pretzel knot complement.
It is a once-punctured torus.  Theorem $7.1$ in \cite{Agol} guarantees that the
punctured torus, which is Fuchsian, remains essential under Dehn filling along the
puncture.  This filling results in an essential torus, which shows that
the resulting manifold cannot be hyperbolike.
\end{proof}

By the Six Theorem, Lemma \ref{upperbound} shows that the longitude must
have length at most six, which combines with the result of Lemma 
\ref{lowerbound} to give that the longitude length for all of these knots
is precisely six.  We obtain the following:

\begin{thm}\label{sixtheoremsharp}
The Six Theorem is sharp for knot complements, with $(p, p, p)$ pretzel
knots all as examples, for every odd $p \geq 3$.
\end{thm}



\end{document}